\documentclass [oneside,a4paper,11pt] {article}
\usepackage{latexsym,epsfig,graphpap,graphics,amsmath,amssymb,}
\begin{document}

\begin{center}
{\bf{COVARIANT ALMOST ANALYTIC VECTOR FIELD ON Q - QUASI UMBILICAL HYPERSURFACE OF A SASAKIAN MANIFOLD}}\\
{\it By}

\vspace{.3cm} \noindent{\bf{Sachin Kumar Srivastava, Alok Kumar Srivastava and Dhruwa Narain}}\\

\end{center}

{\bf Abstract~:}In this paper we have studied the properties of covariant almost analytic vector field on Q - quasi umbilical hypersurface $M$ of a Sasakian manifold $\tilde M$ with $(\phi, g, u, v, \lambda)-$structure and obtained the scalars $\alpha$ and $\beta$ using $1 - form u, v$ covariant almost analytic for the hypersurface $M$ to be totally umbilical and cylindrical.

\vspace{.3cm}\noindent{\bf{2000 Mathematics Subject Classification :}}~ 53D10, 53C25, 53C21.

\vspace{.5cm}
\noindent{\bf{1.~~Introduction~:}}~Let $\tilde M$ be a $(2n + 1)-$dimensional Sasakian manifold with a tensor field $\phi$ of type (1, 1),  a fundamental vector field $\xi$ and $1-form \, \eta$ such that\\
\nolinebreak
(1.1)\hspace{1in}$\eta(\xi) = 1$\\
(1.2)\hspace{1in}$\tilde \phi^2 = - I + \eta \otimes \xi $\\
\vspace{.2cm}
where $I$ denotes the identity transformation.\\
(1.3)(a)\hspace{.8in}
$\eta \, o\,\phi = 0\qquad (b) \qquad \phi \xi = 0 \qquad (c) \qquad rank(\tilde \phi) = 2n$\\
\vspace{.2cm}
If $\tilde M$ admits a Riemannian metric $\tilde g$, such that\\
\vspace{.2cm}
(1.4)\hspace{1in}$\tilde g(\tilde\phi X,\tilde\phi Y) = \tilde g(X, Y) - \eta(X) \eta(Y)$\\
\vspace{.2cm}
(1.5)\hspace{1in}$ \tilde g(X,\xi)=\eta(X)$\\
\vspace{.2cm}
then $\tilde M$ is said to admit a $(\tilde\phi,\xi,\eta,\tilde g)-$ structure called contact metric strucure.\\
If moreover,\\
(1.6)\hspace{1in}$(\tilde\nabla_X \tilde\phi )\, Y =\tilde g (X, \, Y) \xi - \eta (Y) X $\\
 \hspace{.3in} and\\
(1.7)\hspace{1in}$\tilde\nabla_X\xi=-\tilde\phi X$\\
where $\tilde\nabla$ denotes the Riemannian connection of the Riemannian metric g, then $(\tilde M,\tilde\phi,\xi,\eta,\tilde g)$ is called a Sasakian manifold [9].
\vspace{3mm}\parindent=8mm
If we define ${}^\prime F(X, Y) = g (\phi X, Y)$, then in addition to above relation we find\\
(1.8)\hspace{1in}${}^\prime F(X, Y) + {}^\prime F(Y, X) = 0 $\\
(1.9)\hspace{1in}${}^\prime F(X, \phi Y) = {}^\prime F(Y, \phi X)$\\
(1.10)\hspace{.92in}${}^\prime F(\phi X, \phi Y) = {}^\prime F(X, Y)$\\
\vspace{.5cm}
\noindent{\bf{2.~~Hypersurface of a Sasakian manifold with $(\phi, g, u, v, \lambda)-$structure~:}}~\\
\vspace{.1cm}
\hspace{.5cm}Let us consider a $2n-$dimensional manifold $M$ embedded in $\tilde M$ with embedding $b : M \rightarrow \tilde M$. The map $b$ induces a linear transformation map $B$ (called Jacobian map), $B : T_p \rightarrow T_{b_p}$.
\vspace{.1cm}
Let an affine normal $N$ of $M$ is in such a way that $\tilde\phi N$ is always tangent to the hypersurface and satisfying the linear transformations

\vspace{.25cm}\noindent
(2.1) \hspace{1in} $\tilde\phi BX = B \phi X + u(X) N$\\
(2.2) \hspace{1in} $\tilde\phi N = - BU$\\
(2.3) \hspace{1in} $\xi = BV + \lambda N$\\
(2.4) \hspace{1in} $\eta(BX) = v(X)$

\vspace{.5cm}\noindent
where $\phi$ is a (1, 1) type tensor; $U,\, V$ are vector fields; $u, v$ are $1- form$ and $\lambda$ is a $C^\infty -$ function. If $u \ne 0$, $M$, is called a noninvariant hypersurface of $\tilde M$ [1].

\vspace{.3cm}\parindent=8mm
Operating (2.1), (2.2), (2.3) and (2.4)  by $\tilde\phi$ and using (1.1), (1.2) and (1.3) and taking tangent normal parts separately, we get the following induced structure on $M$,

\vspace{.3cm}\noindent
(2.5)(a) \hspace{.96in} $\phi^2 X = -  X + u (X) U + v (X) V$

\vspace{.3cm}\parindent=8mm
(b)\hspace{1in} $u(\phi X) = \lambda v (X), \qquad v(\phi X) =  -\, \eta (N) \, u(X)$

\vspace{.3cm}\parindent=8mm
(c)\hspace{1in} $\phi U = - \, \eta(N) \, V , \qquad \phi V = \lambda U$

\vspace{.3cm}\parindent=8mm
(d)\hspace{1in} $u(U) = 1 -  \lambda \eta(N), \qquad u(V) = 0$

\vspace{.3cm}\parindent=8mm
(e)\hspace{1in} $v(U) = 0, \qquad v(V) = 1 -   \lambda \eta(N)$

\vspace{.3cm}\noindent
and from (1.4) and (1.5), we get the induced metric $g$ on $M$, i.e.,

\vspace{.3cm}
\noindent
(2.6) \hspace{1in} $g(\phi X, \phi Y) = g (X, Y) - u(X) u(Y) - v(X) v(Y)$\\
(2.7) \hspace{1in} $g (U, X) = u (X), \qquad g (V, X) = v (X).$

\vspace{.3cm}
\parindent=8mm
If we consider $\eta(N) = \lambda$, we get the following structures on $M$.

\vspace{.3cm}
\noindent
(2.8)(a)\hspace{1.1in} $\phi^2 = -  I +  u \otimes U + v \otimes V$

\vspace{.3cm}
\parindent=8mm
(b) \hspace{1in} $\phi U = - \, \lambda V, \qquad \phi V = \lambda U$

\vspace{.3cm}
\parindent=8mm
(c) \hspace{1in} $u \, \circ \phi = \lambda v, \qquad v \, \circ \, \phi =  - \lambda u$

\vspace{.3cm}
\parindent=8mm
(d) \hspace{1in} $u( U) = 1 - \lambda^2 , \qquad u(V) = 0$

\vspace{.3cm}
\parindent=8mm
(e) \hspace{1in} $v (U) = 0,  \qquad v (V) = 1 - \lambda^2.$

\vspace{.3cm}
\parindent=8mm
A manifold $M$ with a metric $g$ satisfying (2.6), (2.7) and (2.8) is called manifold with $(\phi, g, u, v, \lambda )-$structure[2].
\vspace{.3cm}
\parindent=8mm
Let $\nabla$ be the induced connection on the hypersurface $M$ of the affine connection $\tilde\nabla$ of $\tilde M$.

\vspace{.3cm}
\parindent=8mm
Now using Gauss and Weingarten's equations

\vspace{.3cm}
\noindent
(2.9)\hspace{1in} $\tilde\nabla _{BX} BY = B\nabla_X Y + h (X, Y) N$

\vspace{.3cm}
\noindent
(2.10) \hspace{.9in} $\tilde\nabla_{BX} N = BHX + w (X) N,$ ~where~ $g (HY, Z) = h(Y, Z).$

\vspace{.3cm}
\parindent=8mm
Here $h$ and $H$ are the second fundamental tensors of type (0, 2) and (1, 1) and $w$ is a $1-form$. Now differentiating (2.1), (2.2), (2.3) and (2.4) covariantly and using (2.9), (2.10), (1.6) and reusing (2.1), (2.2), (2.3) and (2.4), we get

\vspace{.3cm}
\noindent
(2.11) \hspace{.9in} $(\nabla_Y \phi) (X) = v (X) Y - g (X, Y) V - h (X, Y) U - u (X) HY$

\vspace{.3cm}
\noindent
(2.12) \hspace{.9in} $(\nabla_Y u) (X) = - h(\phi X, Y) - u (X) w (Y) - \lambda g (X, Y)$

\vspace{.3cm}
\noindent
(2.13) \hspace{.9in} $(\nabla_Y v) (X) = g(\phi Y, X) + \lambda h (X, Y)$

\vspace{.3cm}
\noindent
(2.14) \hspace{.9in} $\nabla_Y U = w (Y) U - \phi HY - \lambda Y$

\vspace{.3cm}
\noindent
(2.15) \hspace{.9in} $\nabla_Y V = \phi Y + \lambda H Y$

\vspace{.3cm}
\noindent
(2.16) \hspace{.9in} $h(Y, V) = u(Y) - Y \lambda - \lambda w(Y)$

\vspace{.3cm}
\noindent
(2.17) \hspace{.9in} $h(Y, U) = - u(HY)$

\vspace{.3cm}
\parindent=8mm
Since $h (X, Y) = g (HX, Y)$, then from (1.5), and (2.17), we get

\vspace{.3cm}
\noindent
(2.18) \hspace{.9in} $h (Y, U) = 0 \Rightarrow HU = 0.$\\\\
\vspace{.1cm}
\noindent{\bf{3.~~Q - Quasi Umbilical hypersurface~:}}~If

\vspace{.3cm}
\noindent
(3.1) \hspace{1in} $h(X, Y) = \alpha g (X, Y) + \beta q (X) q(Y)$

\vspace{.3cm}
\noindent
where $\alpha, \beta$ are scalar functions, $q$ is $1-form$, then $M$ is called Quasi-umbilical hypersurface and if $g(Q, X) = q(X)$, where $Q$ is vector field, then $M$ is called $Q-$ quasi umbilical hypersurface. If $\alpha = 0$, $\beta \ne 0$, then $Q-$quasi umbilical hypersurface $M$ is called cylindrical hypersurface. If $\alpha \ne 0$, $\beta = 0$, then $Q-$quasi-umbilical hypersurface $M$ is called totally umbilical and if $\alpha = 0, \beta = 0$, then $Q-$quasi umbilical hypersurface is totally geodesic [9].

\vspace{.3cm}
\parindent=8mm
Using (3.1) in (2.11), (2.12), (2.13), (2.14), (2.15), (2.16) and (2.17) we get\\\\
\vspace{.2cm}
(3.2)\hspace{.5in} $(\nabla_Y \phi) (X) = v (X) Y - g (X, Y) V - \{\alpha g (X, Y) + \beta q (X) q(Y)\} U $
$$\hspace*{2.2in} - u (X) \{\alpha Y + \beta q (Y) Q\}$$
\vspace{.2cm}
(3.3) \hspace{.5in} $(\nabla_Y u) (X) = - \{\alpha g (\phi X, Y) + \beta q (\phi X) q(Y)\}  - u (X) w (Y) - \lambda g (X, Y)$\\
\vspace{.2cm}
(3.4) \hspace{.5in} $(\nabla_Y v) (X) = g(\phi Y, X) + \lambda \{\alpha g (X, Y) + \beta q (X) q(Y)\}$\\
\vspace{.2cm}
(3.5) \hspace{.5in} $\nabla_Y U = w (Y) U - \{\alpha \phi Y + \beta q(Y) Q\} - \lambda Y$\\
\vspace{.2cm}
(3.6) \hspace{.5in} $\nabla_Y V = \phi Y + \lambda \{\alpha Y + \beta q(Y) Q\}$\\
\vspace{.2cm}
(3.7) \hspace{.5in} $h(Y, V) = \alpha g (V, Y) + \beta q (V) q(Y)$\\
(3.8) \hspace{.5in} $|u(Q)|^2 = - \cfrac{\alpha}{\beta} (1 - \lambda^2)$\\
Also from (3.1), (2.16) and (2.8)(d) we get
\vspace{.2cm}
$$w(U) = \frac{1 - \lambda^2}{\lambda} - \frac{U \lambda}{\lambda}\leqno (3.9).$$

\vspace{.3cm}
\noindent
{\bf{4.~~Covariant almost analytic vector field on $Q-$quasi umbilical hypersurface~:}}~$1-form$ $u$ and $v$ are said to be covariant almost analytic if

\vspace{.3cm}
\noindent
(4.1) \hspace{1in} $u \{(\nabla_X \phi)(Y) - (\nabla_Y \phi) (X)\} = (\nabla_{\phi X} u )(Y) - (\nabla_X u) (\phi Y)$\\
and

\vspace{.3cm}
\noindent
(4.2) \hspace{1in} $v \{(\nabla_X \phi)(Y) - (\nabla_Y \phi) (X)\} = (\nabla_{\phi X} v )(Y) - (\nabla_X v) (\phi Y)$

\vspace{.3cm}
\noindent
{\bf {Theorem~4.1~:}}~{\it{On $Q-$quasi umbilical hypersurface $M$ with $(\phi, g, u, v, \lambda)-$ structure of a Sasakian manifold $\tilde M$, if $1-form$ $u$ is covariant almost analytic, then we have}}

\vspace{.3cm}
\noindent
(4.3) \hspace{.2in} $v (Y) u(X) - v(X) u(Y) = - 2 \alpha g (\phi X, \phi Y) - 2 \lambda g (\phi X, Y)- \,\beta q (\phi Y) q (\phi X)$\\
$\hspace*{2.3in}  - \, \beta q (Y) q (X) + 2 \beta q (X) u (Y) u (Q)$\\
$\hspace*{2.3in} - \, \beta u (X) q (Y) - u (Y) w (\phi X) + u (\phi Y) w (X)$

\vspace{.3cm}
\noindent
{\bf {Proof~:}}~ From (3.2), we have
$$(\nabla_Y \phi) (X) = v(X) Y - g (X, Y) V - \{\alpha g (X, Y) + \beta q (X) q (Y) \} \, U$$
$$\hspace*{1.8in}- \, u (X) \{\alpha Y + \beta q (Y) Q\}$$
and \hspace{.5in}$(\nabla_X \phi) (Y) = v(Y) X - g (Y, X) V - \{\alpha g (Y, X) + \beta q (Y) q (X) \} \, U$\\
$\hspace*{2.8in}- \, u (Y) \{\alpha X + \beta q (X) Q\}$\\
$.^.. \hspace{.7in} (\nabla_X \phi) (Y) - (\nabla_Y \phi) (X) = v (Y) X - v (X) Y - u (Y) \{\alpha X + \beta q (X) Q\}$ \\
$\hspace*{2.8in} + \, u (X) \{\alpha Y + \beta q (Y) Q\}$

\vspace{.3cm}
\noindent
(4.4)\qquad \quad $u \{(\nabla_X \phi) (Y) - (\nabla_Y \phi) (X)\} = v (Y) u(X) - v (X) u (Y)$\\
$\hspace*{3.2in} + \, \beta \{u (X) q (Y) - q (X) u(Y)\}\, u (Q)$

\vspace{.3cm}
\noindent
Also from (3.4), we have

\vspace{.3cm}
\noindent
\hspace*{.5in} $(\nabla_X u) (Y) = - \{\alpha g (\phi Y, X) + \beta q (\phi Y) q(X)\}  - u (Y) w (X) - \lambda g (X, Y)$

\vspace{.3cm}
\noindent
or \hspace{.3in} $(\nabla_{\phi X} u) (Y) = - \,\alpha g (\phi X, \phi Y) - \beta q (\phi Y) q(\phi X)  - u (Y) w (\phi X) - \lambda g (\phi X, Y)$

\vspace{.3cm}
\noindent
and \hspace{.3in} $(\nabla_X u) (\phi Y) = - \,\alpha g (\phi^2Y, X) - \beta q (\phi^2Y) q(X)  - u (\phi Y) w (X) - \lambda g (X, \phi Y)$

\vspace{.3cm}
\noindent
or \hspace{.3in} $(\nabla_X u) (\phi Y) = \alpha g (Y, X) - \alpha u (X) u(Y) - \alpha v (X) v (Y) + \beta q (Y) q (X)$ \\
\hspace*{1.0in}$-\, \beta u (Y) q (U) q(X) - \beta v (Y) q (V) q (X) - u (\phi Y) w(X) - \lambda g (X, \phi Y)$

\vspace{.3cm}
\noindent
or \hspace{.3in} $(\nabla_X u) (\phi Y) = \alpha g (\phi X, \phi Y) + \beta q (Y) q (X) -  \beta u (Y) q(U) q (X) $ \\
\hspace*{1.2in}$-\, u (\phi Y) w(X) - \lambda g (X, \phi Y)$

\vspace{.3cm}
\noindent
(4.5) \hspace{.1in} $(\nabla_{\phi X} u) (Y) - (\nabla_X u) (\phi Y) = - 2 \alpha g (\phi X, \phi Y) - \beta q (\phi Y) q (\phi X) - \beta q (Y) q (X)$\\
$\hspace*{2.2in}  + \, \beta u (Y) q(U) q (X) - u (Y)w(\phi X) + u (\phi Y) w(X)$\\
$\hspace*{2.2in} - \, \lambda g (\phi X, Y) +  \lambda g (X, \phi Y)$

\vspace{.3cm}
\noindent
using (4.1) in (4.4) and (4.5), we get (4.3).

\vspace{.3cm}
\noindent
{\bf{Corollary~4.1~:}}~{\it{On $Q-$quasi umbilical hypersurface $M$ with $(\phi, g, u, v, \lambda)-$ structure of a Sasakian manifold $\tilde M$, if $1-form$ $u$ is covariant almost analytic, then we  also have}}\\
\vspace{.3cm}
\noindent
(4.6) \hspace{1.5in} $\nabla_U \lambda = 2 \lambda^2.$

\vspace{.3cm}
\noindent
{\bf{Proof~:}}~Put $X = U$ and $Y = V$ in (4.3) and using (2.8)(d), (2.8)(e), we get
$$(1 - \lambda^2)^2 \,= \,- 2 \,\alpha . 0 + 2 \lambda^2 (1 - \lambda^2) - \beta .0 - 0 + \lambda (1 - \lambda^2) \, w (U).$$
Using (3.9), we get

\vspace{.3cm}
\noindent
\hspace*{1.3in} $(1 - \lambda^2) = 2 \lambda^2 + \lambda \Bigg\{\cfrac{1 - \lambda^2}{\lambda} - \cfrac{(U \lambda)}{\lambda}\Bigg\}$

\vspace{.3cm}
\noindent
or \hspace{1.1in} $(1 - \lambda^2) - 2 \, \lambda^2 - 1 + \lambda^2 = - \, U \lambda $

\vspace{.3cm}
\noindent
or \hspace{1.1in} $2\,\lambda^2 = U \lambda$

\vspace{.3cm}
\noindent
or \hspace{1.1in} $U \lambda = 2\,\lambda^2  $

\vspace{.3cm}
\noindent
or \hspace{1.1in} $\nabla_U \lambda = 2 \lambda^2.$

\vspace{.3cm}
\noindent
{\bf{Theorem~4.2~:}}~{\it {On cylindrical hypersurface $M$ with $(\phi, g, u, v, \lambda)-$structure of a Sasakian manifold $\tilde M$, with covariant almost analytic vector fields $U$, $V$, we have}}

\vspace{.3cm}
\noindent
(4.7) \hspace{.3in} $v(Y) u(X) - v(X) u (Y) = - 2 \lambda g (\phi X, Y) - \beta q (\phi X) q (\phi Y) - \beta q (Y) q (X)$\\
$\hspace*{2.5in}  - \, \beta u (X) q (Y) - u (Y)w(\phi X) + u (\phi Y) w(X)$

\vspace{.3cm}
\noindent
{\bf{Proof~:}}~Putting $\alpha = 0$ in (4.3) and using $u(Q) = 0$, we get (4.7).

\vspace{.3cm}
\noindent
{\bf{Corollary~4.2~:}}~{\it {On cylindrical hypersurface $M$ with $(\phi, g, u, v, \lambda)-$structure of a Sasakian manifold $\tilde M$, with covariant almost analytic vector fields $U$, $V$, we also have}}

\vspace{.3cm}
\noindent
(4.8) \hspace{1in} $V \lambda = 0$, i.e. $\lambda$ is covariant constant along $V$

\vspace{.3cm}
\noindent
(4.9) \hspace{1.5in} $q (Q) = 0$.

\vspace{.3cm}
\noindent
{\bf{Proof~:}}~Putting $X = V$ in (4.7), we get

\vspace{.3cm}
\noindent
\hspace*{.5in} $-\,(1 - \lambda^2) u(Y) = - \,2 \lambda g (\phi V, Y) - \beta q (\phi Y) q (\phi V) - \beta q (V) q(Y)$\\
\hspace*{2.3in} $-\,\beta u (V) q(Y) - u(Y) w (\phi V) + u (\phi Y) w (V)$

\vspace{.3cm}
\noindent
or \hspace{.3in} $-\,(1 - \lambda^2) u(Y) = - \,2 \lambda^2 u (Y) - \lambda u(Y) w(U) + u(\phi Y) w(V) $

\vspace{.3cm}
\noindent
or \hspace{.3in} $(3 \lambda^2 -\,1) u (Y) = -\, \lambda u(Y) \Bigg(\cfrac{1 - \lambda^2}{\lambda} - \cfrac{U \lambda}{\lambda}\Bigg) + u(\phi Y)\Bigg(\cfrac{ - V\,\lambda}{\lambda}\Bigg) $

\vspace{.3cm}
\noindent
or \hspace{.3in} $(3 \lambda^2 -\,1) u (Y) = -\, u(Y) \{(1 - \lambda^2) - (U \lambda)\} - u(\phi Y)\Bigg(\cfrac{  V\,\lambda}{\lambda}\Bigg)$

\vspace{.3cm}
\noindent
or \hspace{.3in} $(3 \lambda^2 -\,1 + 1 - \lambda^2) U =  (U \lambda) U + \Bigg(\cfrac{  V\,\lambda}{\lambda}\Bigg)\, \phi U$

\vspace{.3cm}
\noindent
or \hspace{.3in} $2 \lambda^2\, U =  (U \lambda)\, U + \Bigg(\cfrac{  V\,\lambda}{\lambda}\Bigg)\, \phi U$

\vspace{.3cm}
\noindent
using (4.6), we get $V \lambda = 0$, i.e. $\lambda $ is covariant constant along $V$.

\vspace{.3cm}
\noindent
Also from (4.7), we get
$$u(X) V - v(X) U = - \, 2 \lambda \phi X + \beta \phi Q q (\phi X) - \beta q (X) Q - \beta u (X) Q $$
$$\hspace*{1.5in} -\, w(\phi X) \, U - w\, (X) \phi U$$

\vspace{.3cm}
\noindent
Contracting with respect to $X$, we get
$$0 = \beta q (\phi^2 Q) - \beta q(Q) - \beta u(Q) + \lambda w (V) + \lambda w (V)$$
or \hspace{1.5in} $\beta q (Q) = - (V \lambda) = 0$

\vspace{.3cm}
\noindent
since $\beta \ne 0$, so $q(Q) = 0.$

\vspace{.3cm}
\noindent
{\bf{Theorem~4.3~:}}~{\it {On $Q-$quasi umbilical hypersurface $M$ with $(\phi, g, u, v, \lambda)-$structure of a Sasakian manifold $\tilde M$, where $ 1 - form u$ is covariant almost analytic. If $M$ is totally umbilical, then we have}}

\vspace{.3cm}
\noindent
(4.10) \hspace{.5in} $v(Y) u(X) - v(X) u(Y) = - \, 2 \alpha g (\phi X, \phi Y) - 2 \lambda g (\phi X, Y)$\\
\hspace*{3.2in} $-\, u (Y) w (\phi X) + u (\phi Y) w(X)$

\vspace{.3cm}
\noindent
(4.11) \hspace{.5in} $\alpha = - \, \frac{(V \lambda)}{1 + \lambda^2}$ \quad or \quad $\alpha = - \, d (\tan^{-1} \lambda) \, V$.

\vspace{.3cm}
\noindent
{\bf{Proof~:}}~Putting $\beta = 0$ in (4.3), we get (4.10). Further putting $X = U$ in (4.10), we get

\vspace{.3cm}
\noindent
\hspace*{.2in} $(1 - \lambda^2) v(Y) = 2\, \alpha \lambda v (\phi Y) + 2 \lambda^2 v (Y) + \lambda u (Y) w (V) + u (\phi Y) w(U)$

\vspace{.3cm}
\noindent
\hspace*{.2in} $(1 - \lambda^2) v(Y) = - \,2\, \alpha \lambda^2 u (Y) + 2 \lambda^2 v(Y) $\\
\hspace*{1.5in} $+\,\lambda u(Y) \Bigg\{\alpha \cfrac{(\lambda^2 - 1 )}{\lambda} - \cfrac{(V \lambda )}{\lambda}\Bigg\} + u(\phi Y)\Bigg\{\cfrac{1 - \lambda^2}{\lambda} - \cfrac{(U \lambda)}{\lambda}\Bigg\} $

\vspace{.3cm}
\noindent
or \hspace{.3in} $(1 - \lambda^2 -\,2 \lambda^2 - 1 + \lambda^2) v (Y) = \{-\, 2 \alpha \lambda^2 + \alpha \lambda^2 - \alpha \}\,u (Y) - u(Y) (V \lambda)$
\hspace*{3.3in}$- v(Y) (U \lambda) $

\vspace{.3cm}
\noindent
or \hspace{.3in} $-\,2 \lambda^2 V  = -\, \alpha \,(1 + \lambda^2)U - (V \lambda)\, U -  (U\,\lambda )\, V$

\vspace{.3cm}
\noindent
Operating by $u$

\vspace{.3cm}
\noindent
or \hspace{1in} $0 = -\, \alpha\,(1 + \lambda^2)\,(1 - \lambda^2) - (V \lambda)\, (1 - \lambda^2)$

\vspace{.3cm}
\noindent
or \hspace{1in} $\alpha = - \cfrac{(V\,\lambda)}{1 + \lambda^2}$

\vspace{.3cm}
\noindent
or \hspace{1in} $\alpha = - \, d\,(\tan^{-1} \lambda) \, (V)$

\vspace{.3cm}
\noindent
{\bf{Corollary~4.3~:}}~{\it{On $Q-$quasi umbilical hypersurface $M$ with $(\phi, g, u, v, \lambda)-$structure of a Sasakian manifold $\tilde M$, where $ 1 - form $ $u$ is covariant almost analytic. If $M$ is totally umbilical, then we also have}}

\vspace{.3cm}
\noindent
(4.12) \hspace{.5in} $\alpha = - \,\cfrac{1}{\sqrt{2n - 1}}\, d \,\Bigg(\tan^{- 1} \, \cfrac{\lambda} {\sqrt{2n - 1}}\Bigg) \, (V)$

\vspace{.3cm}
\noindent
{\bf{Proof~:}}~From (4.10), we have
$$u(X) V - v(X) U = 2 \alpha \phi^2 X - 2 \lambda \phi X - w (\phi X) U + \lambda w (X) \, V$$
Contracting with respect to X

\vspace{.3cm}
\noindent
\hspace*{.5in} $0 = 2\, \alpha \, \{- 2n + 4 \alpha \,(1 - \lambda^2)\} +  2 \lambda w (V)$

\vspace{.3cm}
\noindent
or \hspace{.3in} $0 = - \,2\, \alpha \,n + 2 \alpha \,(1 - \lambda^2) + \lambda w (V) $

\vspace{.3cm}
\noindent
or \hspace{.3in} $- \, 2\, \alpha \,n + 2 \alpha \,(1 - \lambda^2) + \lambda \Bigg\{\alpha \cfrac{(\lambda^2 - 1) }{\lambda} - \cfrac{(V \lambda)}{\lambda}\Bigg\}  = 0$

\vspace{.3cm}
\noindent
or \hspace{.3in} $\alpha  \,(1 - \lambda^2 - 2n ) = (V \lambda)$

\vspace{.3cm}
\noindent
or \hspace{.3in} $\alpha = \cfrac{V \lambda}{1 - \lambda^2 - 2n } = -\, \cfrac{(d \lambda) \, V}{(2n - 1) + \lambda^2}$

\vspace{.3cm}
\noindent
thus \hspace{.15in} $\alpha = - \,\cfrac{1}{\sqrt{2n - 1}}\, d \,\Bigg(\tan^{- 1} \, \cfrac{\lambda} {\sqrt{2n - 1}}\Bigg) \, (V).$

\vspace{.3cm}
\noindent
{\bf{Theorem~4.4~:}}~{\it{On a $Q-$quasi umbilical noninvariant hypersurface $M$ with $(\phi, g, u,$ $ v, \lambda)-$structure of a Sasakian manifold $\tilde M$, If $1-form$ $v$ is covariant almost analytic, then we have}}

\vspace{.3cm}
\noindent
(4.13) \hspace{.5in} $\alpha \{u(X) v(Y) - u(Y) v(X)\} = - \, 2 g (\phi X, \phi Y) + 2 \lambda \alpha g (\phi X, Y) $\\
\hspace*{3.2in} $+\,\lambda \beta \{q(\phi X) q (Y) - q (X) q (\phi Y) \}$

\vspace{.3cm}
\noindent
{\bf{Proof~:}}~From (3.2), we get
$$(\nabla_X \phi) (Y) - (\nabla_Y \phi) (X) = v (Y) X - v(X) Y  - \, u(Y)\{\alpha X + \beta q (X) Q\} $$
$$\hspace*{1.3in} +\, u (X) \{\alpha Y + \beta q (Y) Q\}$$

\noindent
(4.14) \hspace{.3in} $v\, \{(\nabla_X \phi) (Y) - (\nabla_Y \phi) (X)\} = \, \alpha \{u (X)\,v(Y) - u (Y) \, v(X)\}$

\vspace{.3cm}
\noindent
from (3.4), we get
$$(\nabla_{\phi X} v) (Y) = -\, g(X, Y) + u(X) u(Y) + v (X) v(Y) + \lambda \alpha g (\phi X, Y) + \lambda \beta q (\phi X) q(Y)$$
and
$$(\nabla_X v) (\phi Y) = g(X, Y) - u(X) u(Y) - v (X) v(Y) + 2\lambda g (X, \phi Y) + \lambda \beta q (X) q(\phi Y)$$
therefore

\vspace{.3cm}
\noindent
(4.15) \hspace{.1in} $(\nabla_{\phi X} v) (Y) - (\nabla_X v) (\phi Y) = -\,2\,g(\phi X, \phi Y) + \lambda \beta \{q (\phi X) q(Y) - q(X) q(\phi Y)\} $\\
\hspace*{3.3in} $+ 2 \lambda \alpha g (\phi X, Y)$

\vspace{.3cm}
\noindent
Using (4.2) in (4.14) and (4.15), we get (4.13).

\vspace{.3cm}
\noindent
{\bf{Theorem~4.5~:}}~{\it{Let $1-form$ $v$ is covariant almost analytic on $Q-$quasi umbilical noninvariant hypersurface $M$ with $(\phi, g, u, v, \lambda)-$structure of a Sasakian manifold $\tilde M$, and if it is cylindrical also, we have}}

\vspace{.3cm}
\noindent
(4.16) \hspace{.5in} $2 g (\phi X, \phi Y) = \,\lambda \beta \{q(\phi X) q (Y) - q (X) q (\phi Y) \}$, with $X \ne U$.

\vspace{.3cm}
\noindent
{\bf{Proof~:}}~Putting $\alpha = 0$ in (4.13), we get (4.16). Further if $X = U$, then from (4.16), we get
$$2 g (\phi U, \phi Y) = \,\lambda \beta \{- \, \lambda q(V) q (Y) - q (U) q (\phi Y) \}$$
\noindent
or \hspace{1.05in} $- \, 2 \lambda g \, (V, \phi Y) = 0$

\vspace{.3cm}
\noindent
or \hspace{1.05in} $- \, 2 \lambda v \, (\phi Y) = 0$

\vspace{.3cm}
\noindent
since $v(\phi Y) \ne 0$ and $\lambda \ne 0$, therefore $X \ne U$.

\vspace{.3cm}
\noindent
{\bf{Theorem~4.6~:}}~{\it{Let $1-form$ $v$ is covariant almost analytic on $Q-$quasi umbilical noninvariant hypersurface $M$ with $(\phi, g, u, v, \lambda)-$structure of a Sasakian manifold, and if it is totally umbilical, then we have}}

\vspace{.3cm}
\noindent
(4.17) \hspace{.5in} $\alpha \{u(X) v(Y) - v(X) u(Y) \} = -2\,g (\phi X, \phi Y) + 2 \lambda \alpha g (\phi X, Y)$, with $X \ne U$.

\vspace{.3cm}
\noindent
{\bf{Proof~:}}~Putting $\beta = 0$ in (3.13), we get (3.17). If we take $X = U$ in (4.17), we get
$$\alpha \, (1 - \lambda^2)\, v(Y) = 2\,\lambda v (\phi Y) - 2\, \lambda \alpha v (Y)$$
\noindent
or \hspace{1.3in} $\alpha \, \{1 - \lambda^2 + 2 \lambda\}\, v(Y) = 2\,\lambda v (\phi Y)$

\vspace{.3cm}
\noindent
or \hspace{1.3in} $\alpha \, \{\lambda^2 - 2 \lambda - 1\}\, v(Y) = -\,2\,\lambda \phi V = -\,2 \lambda^2 U$

\vspace{.3cm}
\noindent
i.e. $V$ and $U$ are linearly dependent, which is a contradiction. Thus $X \ne U$.

\begin{center}
{\bf References}
\end{center}

\begin{enumerate}
\item[{[1]}] S. I. Goldberg and K. Yano : {\em Noninvariant hypersurfaces of almost contact manifolds}, J. Math. Soc., Japan, {\bf 22} (1970), 25-34. 
\item[{[2]}] D. Narain : {\em Hypersurfaces with $(f, g, u, v, \lambda)-$structure of an affinely cosymplectic manifold}, Indian Jour. Pure and Appl. Math.{\bf 20} No. 8 (1989), 799-803.
\item[{[3]}] Dhruwa Narain and S. K. Srivastva : {\em Noninvariant Hypersurfaces of Sasakian Space Forms}, Int. J. Contemp. Math. Sciences, Vol. {\bf 4}, 2009, no. 33, 1611 - 1617.
\item[{[4]}] Dhruwa Narain, S. K. Srivastava and Khushbu Srivastava : {\em A note of noninvariant  hypersurfaces of para Sasakian manifold}, IOSR Journal of Engineering(IOSRJEN) Vol.{\bf 2}, issue 2 (2012), 363-368.
\item[{[5]}] Dhruwa Narain and S. K. Srivastava : {\em On the hypersurfaces of almost r-contact manifold}, J. T. S. Vol.{\bf 2}(2008), 67-74.
\item[{[6]}] Rajendra Prasad and M. M. Tripathi : {\em Transversal Hypersurfaces of Kenmotsu Manifold}, Indian Jour. Pure and Appl. Math.{\bf 34} No. 3 (2003), 443-452.
\item[{[7]}] B. B. Sinha and Dhruwa Narain : {\em Hypersurfaces of nearly Sasakian manifold}, Ann. Fac. Sci. Dekishasa Zaire, {\bf 3} (2), (1977), 267-280.
\item[{[8]}] M. Okumura : {\em Totally umbilical hypersurface of a locally product Riemannian manifold}, Kodai Math. Sem. Rep.{\bf 19} (1967), 35-42.
\item[{[9]}] D. E. Blair  : {\em Almost contact manifolds with killing structure tensor}, Pacific J. of Math. {\bf 39} (12), (1971), 285-292.
\item[{[10]}] T. Miyazawa and S. Yamaguchi : {\em Some theorems on K-contact metric manifolds and Sasakian manifolds}, T. R. V. Math. {\bf 2} (1966), 40-52.
\end{enumerate}
\vspace{.2cm}
\hspace{.5cm}{\it{ Authors' addresses}:}\\
{Sachin Kumar Srivastava}\\
{Amity Institute of Applied Sciences,Amity University , Noida, U.P., India}\\
\vspace{.2cm}
$E-mail: sksrivastava4@amity.edu$\\
Alok Kumar Srivastava\\
Department of Mathematics ,Govt. Degree College, Chunar -
Mirzapur,U.P., India \\
\vspace{.2cm}
$E-mail: aalok\_sri@yahoo.co.in$\\
Dhruwa Narain\\
Department of Mathematics and Statistics\\
D.D.U. Gorakhpur University Gorakhpur, Gorakhpur-273009, India\\
$E- mail:dhruwanarain\_dubey@yahoo.co.in$\\

\end{document}